\documentclass[11pt]{article}
\usepackage{amsmath}
\usepackage{amsfonts}
\usepackage{bm}
\usepackage{graphicx}
\usepackage{esdiff}
\usepackage{url}
\usepackage[linesnumbered,ruled,vlined]{algorithm2e}

\usepackage{mathptmx}      % use Times fonts if available on your TeX system
\usepackage{epstopdf}

\newcommand{\R}{{\mathbb R}}

\newcommand{\be}{\begin{equation}}
\newcommand{\ee}{\end{equation}}
\newcommand{\ba}{\begin{array}}
\newcommand{\ea}{\end{array}}
\newcommand{\baa}{\left[\begin{array}}
\newcommand{\eaa}{\end{array}\right]}

\newcommand{\beqa}{\begin{eqnarray}}
\newcommand{\eeqa}{\end{eqnarray}}
\newcommand{\bt}{\begin{tabular}}
\newcommand{\et}{\end{tabular}}

\newcommand{\bi}{\begin{itemize}}
\newcommand{\ei}{\end{itemize}}
\newcommand{\bc}{\begin{center}}
\newcommand{\ec}{\end{center}}

\newtheorem{prop}{Proposition}
\newtheorem{remark}{Remark}

\def\QED{\hfill \mbox{\rule[0pt]{1.5ex}{1.5ex}}}
\title{On the Relation Between the Randomized Extended Kaczmarz Algorithm
and Coordinate Descent}

\author{Bogdan Dumitrescu$^{\rm a,b}$}
\begin{document}
\maketitle

\vspace*{20mm}
\begin{abstract}
In this note we compare the randomized extended Kaczmarz (EK) algorithm and
randomized coordinate descent (CD) for solving the full-rank
overdetermined linear least-squares problem and prove that
%CD is always faster than EK.
CD needs less operations for satisfying the same residual-related termination criteria.
For the general least-squares problems, we show that running first CD to
compute the residual and then standard Kaczmarz
on the resulting consistent system is more efficient than EK.

\vspace*{5mm}
{\bf Keywords:} randomized algorithms, least-squares, Kaczmarz method, coordinate descent

\vspace*{5mm}
{\bf MSC:} 65F10, 65F20, 15A06
\end{abstract}

\vspace*{50mm}
\noindent $^{\rm a}$Department of Automatic Control and Computers,
University Politehnica of Bucharest, Spl.\ Independen\c{t}ei 313, Bucharest 060042, Romania.
E-mail: bogdan.dumitrescu@acse.pub.ro \\
$^{\rm b}$Department of Signal Processing, Tampere University of Technology, Finland.
E-mail: bogdan.dumitrescu@tut.fi

%---------------------------------------------------------------------

%---------------------------------------------------------------------
\pagebreak
\section{Introduction}

Given matrix $\bm{A} \in \R^{m \times n}$ and vector $\bm{b} \in \R^m$,
the linear least-squares (LS) problem consists of finding $\bm{x} \in \R^n$
such that $\| \bm{b} - \bm{A} \bm{x} \|_2$ is minimum.
Unless explicitly stated, we consider the full-rank overdetermined problem,
i.e.\ $m \ge n$ and $\mbox{rank} \bm{A} = n$.
Besides standard solutions based on orthogonal triangularization
or the normal equations, significant recent interest was focused
on randomized algorithms, showing clear benefits for certain categories
of problems, especially for large dimensions and sparse matrices.

There are two main classes of randomized algorithms for the LS problem,
both based on simple iterated projection operations.
In coordinate descent (CD) \cite{LeLe10}, at iteration $k$,
the current residual is projected onto a random column of the matrix $\bm{A}$,
in order to obtain the optimal LS update of a single element of
the current approximation of the solution $\bm{x}^{(k)}$.
In the Kaczmarz algorithm \cite{StVe09}, the solution $\bm{x}^{(k)}$
is projected onto the hyperplane defined by a random row of the matrix $\bm{A}$
and the respective element of $\bm{b}$, thus obtaining the next approximation
$\bm{x}^{(k+1)}$.
Unlike CD, randomized Kaczmarz converges only when the system
$\bm{A}\bm{x} = \bm{b}$ is consistent.
Otherwise, it hovers around the LS solution, within guaranteed bounds \cite{Nee10}.
This behavior is natural, since at each iteration the approximated solution
satisfies exactly an equation of the system $\bm{A}\bm{x} = \bm{b}$,
which is not the case in general for the LS solution.
Convergence to the LS solution can be achieved if sub-optimal steps are used,
see e.g.\ \cite{HaNi90} and the references therein,
and the step length goes to zero;
however, convergence speed may become very slow.

To fix this drawback, the extended Kaczmarz (EK) algorithm,
in both randomized \cite{ZoFr13} and original deterministic \cite{Popa95} forms,
simultaneously builds an approximation of the residual,
such that a consistent system is asymptotically obtained,
and applies Kaczmarz iterations for the current approximation of this system.
Thus, the algorithm converges to the LS solution.

We show in this note that EK consists in fact of CD and Kaczmarz iterations,
thus combining both classes of randomized algorithms.
Since CD can find on its own the LS solution, we argue that EK can
never be faster than CD, neither in terms of number of iterations,
nor in terms of number of operations.
So, we conclude that randomized CD should be preferred over EK
in all situations, for overdetermined LS problems.
We discuss also the general LS problem (not full rank) and show that
EK can be safely replaced by CD followed by the usual Kaczmarz
%for faster convergence.
for better practical behavior.
Other combinations of the algorithms are possible for providing early
estimates of the solution, like EK.

The notation resembles that from \cite{ZoFr13}.
We denote by $\bm{A}^{(i)}$ and $\bm{A}_{(j)}$ the $i$-th row
and $j$-th column of matrix $\bm{A}$, respectively, both seen as column vectors.
The scalar product of two vectors is denoted $\langle \cdot,\cdot \rangle$
and $[m] = \{ 1, \ldots, m\}$.
The $i$-th unit vector is $\bm{e}_i$.
To distinguish between algorithms, we add the subscript
EK, CD or K (the latter for standard Kaczmarz) to variables having
the same meaning, but different values for the three algorithms.
We denote $\bm{x}_o$ the solution of the LS problem
and ${\cal R}(\bm{A})$ the range of $\bm{A}$.
The 2-norm is used for vectors and matrices unless otherwise stated.

%---------------------------------------------------------------------
\section{Extended Kaczmarz vs coordinate descent}

%---------------------------
\begin{algorithm}[t]
\label{algo:rek}
%\SetKwComment{Comment}{}{}

\BlankLine
\KwData{Matrix $\bm{A} \in \R^{m \times n}$, vector $\bm{b} \in \R^m$,
        stopping tolerances $\varepsilon_{CD}$, $\varepsilon_K$
}
\KwResult{Least-squares solution $\bm{x} \in \R^n$ minimizing
$\| \bm{b} - \bm{A} \bm{x} \|$
}
\BlankLine

Initialize $\bm{r}_{CD}^{(0)} = \bm{b}$,
           $\bm{x}_{CD}^{(0)} = 0$,
           $\bm{x}_{EK}^{(0)} = 0$. \\
\For{$k = 0, 1, 2, \ldots$}{
   Pick $j_k \in [n]$ with probability $p_j = \|\bm{A}_{(j)}\|_2^2 / \|\bm{A}\|_F^2$,
    $j \in [n]$ \\
   Find optimal coordinate descent step:
    $\mu = \frac{\langle \bm{r}_{CD}^{(k)},\bm{A}_{(j_k)}\rangle}
                                         {\|\bm{A}_{(j_k)}\|_2^2}$ \\
   Update CD residual: $\bm{r}_{CD}^{(k+1)} = \bm{r}_{CD}^{(k)} - \mu \bm{A}_{(j_k)}$ \\
   Update CD solution: $\bm{x}_{CD}^{(k+1)} = \bm{x}_{CD}^{(k)} + \mu \bm{e}_{j_k}$ \\
   Pick $i_k \in [m]$ with probability $q_i = \|\bm{A}^{(i)}\|_2^2 / \|\bm{A}\|_F^2$,
    $i \in [m]$ \\
   Update EK solution: $\bm{x}_{EK}^{(k+1)} = \bm{x}_{EK}^{(k)} +
     \frac{\langle \bm{b} - \bm{r}_{CD}^{(k+1)}, \bm{e}_{i_k} \rangle
              - \langle \bm{x}_{EK}^{(k)},\bm{A}^{(i_k)}\rangle}
          {\|\bm{A}^{(i_k)}\|_2^2} \bm{A}^{(i_k)}$ \\
   Check every $8 \min(m,n)$ iterations and terminate if both following conditions hold
   \be
   \frac{\| \bm{A}^T \bm{r}_{CD}^{(k)}\|_2}
        {\|\bm{A}\|_F^2 \|\bm{x}^{(k)}\|_2} \le \varepsilon_{CD} 
   \label{stopCD}
   \ee
   \be   
   \frac{\|\bm{b} - \bm{r}_{CD}^{(k)} - \bm{A}\bm{x}^{(k)}\|_2}
        {\|\bm{A}\|_F \|\bm{x}^{(k)}\|_2} \le \varepsilon_K
   \label{stopK}
   \ee
}
\caption{Randomized Extended Kaczmarz}
\end{algorithm}
%---------------------------

Algorithm \ref{algo:rek} shows a slightly modified version of
the randomized EK from \cite{ZoFr13}.
Besides non-significant permutations of the steps and some different explanations,
only step 6 is new here and does not affect the final outcome.
Let us first discuss the structure of the algorithm and explain its
relation with CD.
For further reference, we denote
\be
\bm{r}^{(k)} = \bm{b} - \bm{A} \bm{x}^{(k)}
\label{resid}
\ee
the residual at iteration $k$.

The EK algorithm, as presented in \cite{ZoFr13}, has two intertwined parts.
In the first, a residual is built, converging to the optimal residual
$\bm{b} - \bm{A} \bm{x}_o$ of the LS problem.
At each iteration, a column $j_k$ is picked randomly as in step 3 of
Algorithm 1, and the residual is projected onto the orthogonal complement
of this column, thus obtaining a new residual (smaller in size than the
previous because of the projection).
However, this is exactly what CD does and that is why we denote this residual
$\bm{r}_{CD}^{(k)}$.
Indeed, in CD, the residual is projected onto column $j_k$ in order to
find the optimal update of the $j_k$-th element of the solution, as in step 4
(this projection maximizes the decrease of $\| \bm{b} - \bm{A} \bm{x}^{(k+1)} \|$
if only the $j_k$-th coordinate of $\bm{x}^{(k)}$ is modified).
After updating the solution as in step 6, the new residual from step 5
is indeed
\[
\bm{r}_{CD}^{(k+1)} = \bm{b} - \bm{A}\bm{x}_{CD}^{(k+1)}
= \bm{b} - \bm{A}\bm{x}_{CD}^{(k)} - \mu \bm{A}_{(j_k)}
= \bm{r}_{CD}^{(k)} - \mu \bm{A}_{(j_k)}
\]
and is orthogonal on column $j_k$, as one can easily check by plugging in the
expression of the optimal update $\mu$:
\[
\langle \bm{r}_{CD}^{(k+1)}, \bm{A}_{(j_k)} \rangle =
\langle \bm{r}_{CD}^{(k)} - \mu \bm{A}_{(j_k)}, \bm{A}_{(j_k)} \rangle = 0.
\]
So, steps 3-5 of EK compute the CD residual.
Only one more arithmetic operation per iteration is necessary to update the CD solution,
as in step 6.

We conclude that Algorithm \ref{algo:rek} without steps 7 and 8
is actually the randomized CD algorithm,
which converges to the LS solution, i.e.\ $\bm{x}_{CD}^{(k)} \rightarrow \bm{x}_o$,
$\bm{r}_{CD}^{(k)} \rightarrow \bm{b} - \bm{A} \bm{x}_o$, see \cite{LeLe10},
or \cite{Nest12rcd} for a more general treatment.
(The probabilities $p_j$ are also taken like in the randomized CD.)
In step 9, CD needs only the stopping criterion (\ref{stopCD}),
which shows that the residual has become nearly orthogonal on ${\cal R}(\bm{A})$.
Note that the stopping criterion (\ref{stopK}) is irrelevant for CD, since
$\bm{b} - \bm{r}_{CD}^{(k)} - \bm{A}\bm{x}_{CD}^{(k)} = 0$ by definition (\ref{resid}).
In what follows, we understand by CD the algorithm described in
this paragraph, with $\bm{x}_{CD}^{(k)}$ used in (\ref{stopCD}).

The second part of Algorithm \ref{algo:rek}, steps 7 and 8,
implements a Kaczmarz iteration for
the LS problem
\be
\bm{A} \bm{x} = \bm{b} - \bm{r}_{CD}^{(k)}.
\label{Ksys}
\ee
Since CD converges to the LS solution, the above system becomes asymptotically
consistent and hence the Kaczmarz iterations converge also to the true solution,
as shown in \cite{ZoFr13}.
Both stopping criteria from step 9 are now necessary;
as above, the criterion (\ref{stopCD}) shows that the residual has converged;
the criterion (\ref{stopK}) shows that the Kaczmarz iterations have
converged and hence a solution to (\ref{Ksys}) has been found.
In \cite{ZoFr13}, the tolerances $\varepsilon_{CD}$ and $\varepsilon_K$ are
equal, but we take them different for the sake of generality.
Formally, EK is Algorithm \ref{algo:rek} with $\bm{x}_{EK}^{(k)}$
used in the stopping criteria (\ref{stopCD}) and (\ref{stopK}) of step 9.

\begin{remark} \rm
\label{rem:CDvsEK}
The above presentation of the CD and EK algorithms allows a quick assessment.
CD computes an approximation of the optimal residual of the LS problem
and produces an approximation of the LS solution which always satisfies (\ref{Ksys}).
EK computes an approximation of the LS solution by approximating the solution
of (\ref{Ksys}), a system depending on the CD residual.
So, EK builds an approximation based on the CD approximation;
EK aims to a target that is driven by CD.
By its very principle, EK should need more iterations than CD to 
meet the same stopping criterion.
We give below a formal proof of this fact.
\QED
\end{remark}

\begin{prop} \rm
\label{prop:term}
In average (over the random draw of columns and rows),
CD terminates in less iterations than EK.
Also, CD needs less arithmetic operations than EK.
\end{prop}

{\em Proof.}
It is enough to prove the proposition for a literal implementation of Algorithm 1,
where the random columns $j_k$ are the same for EK and CD.
Then, by averaging, the same relations hold.

We can safely assume that the EK and CD solution approximations
have similar magnitudes, i.e.\ the values $\|\bm{x}_{EK}^{(k)}\|$
and $\|\bm{x}_{CD}^{(k)}\|$ do not make the stopping criterion (\ref{stopCD})
significantly different for EK and CD, especially near convergence,
the LS solution being unique.
The stopping criterion (\ref{stopCD}) depends only on the CD residual,
so EK cannot stop earlier than CD.
As mentioned above, the stopping criterion (\ref{stopK}) is always met for CD
because (\ref{Ksys}) holds exactly for $\bm{x}_{CD}^{(k)}$.
So, again, EK cannot stop earlier than CD.
Hence, CD needs at most the same number of iterations as EK to terminate.

The number of arithmetic operations per iteration is also in favor of CD,
since CD needs only a subset of the operations of EK (step 6 is negligible
with respect to the others).
CD needs about $4m$ operations per iteration (steps 4 and 5 dictate the complexity),
while EK needs about $4m+4n$ (steps 4, 5 and 8).
\QED

\smallskip
Of course, when the CD and EK algorithms are separately implemented,
then the random columns are different and it may happen that EK terminates
faster than CD, due to a more favorable draw of columns.

Proposition \ref{prop:term} describes the relation between CD and EK
from a strictly computational viewpoint, that of algorithm termination.
However, the relation between their residuals can be more precisely qualified.

\begin{prop} \rm
\label{prop:res}
Asymptotically, the residuals of CD and EK satisfy
$\| \bm{r}_{EK}^{(k)} \| \ge \| \bm{r}_{CD}^{(k)} \|$.
\end{prop}

{\em Proof.}
Define
\be
\hat{\bm{r}}_{EK}^{(k)} = \bm{b} - \bm{r}_{CD}^{(k)} - \bm{A}\bm{x}_{EK}^{(k)}
\stackrel{(\ref{resid})}{=} \bm{A} (\bm{x}_{CD}^{(k)} - \bm{x}_{EK}^{(k)})
\label{rek_fake}
\ee
the residual of the system (\ref{Ksys}) that EK actually attempts to solve at iteration $k$.
It results that
\be
\bm{r}_{EK}^{(k)} = \bm{r}_{CD}^{(k)} + \hat{\bm{r}}_{EK}^{(k)}.
\label{rek}
\ee
%Using (\ref{resid}) for CD, it also follows from (\ref{rek_fake}) that
%\[
%\hat{\bm{r}}_{EK}^{(k)} = \bm{A} (\bm{x}_{CD}^{(k)} - \bm{x}_{EK}^{(k)}).
%\]
Since CD converges to the solution of the LS problem,
its residual tends to become orthogonal to the range of $\bm{A}$,
thus (\ref{rek_fake}) implies that
\[
\langle \bm{r}_{CD}^{(k)}, \hat{\bm{r}}_{EK}^{(k)}\rangle \rightarrow 0.
\]
Hence one can infer from (\ref{rek}) that, asymptotically,
$\| \bm{r}_{EK}^{(k)} \| \ge \| \bm{r}_{CD}^{(k)} \|.$
\QED

\smallskip
The above Propositions show that CD reaches its target faster than EK.
Of course, a smaller residual does not necessarily mean that the solution
approximation is closer to the LS optimum, although this is more likely.
In this context, one may wonder if the convergence speed is indeed
different for the two algorithms.

\begin{remark} \rm
\label{rem:convCD-EK}
In \cite[Th.2.3]{ZoFr13}, the CD residual is shown to satisfy the relation
\be
\label{convCD}
E\{ \| \bm{r}_{CD}^{(k)} - \bm{r}_o \|^2 \} \le
   \left( 1 - \frac{1}{\kappa_F^2( \bm{A} )} \right)^k \| \bm{b} - \bm{r}_o \|^2,
\ee
where $\bm{r}_o = \bm{b} - \bm{A} \bm{x}_o$ is the optimal residual and
$\kappa_F(\bm{A}) = \| \bm{A} \|_F \| \bm{A}^\dagger \|$, with $\bm{A}^\dagger$
the Moore-Penrose pseudoinverse of $\bm{A}$.
The average in (\ref{convCD}) is taken over the random column indices generated
in step 3 of Algorithm \ref{algo:rek}.
Since
\[
\bm{A} (\bm{x}_{CD}^{(k)} - \bm{x}_o ) = \bm{r}_o - \bm{r}_{CD}^{(k)},
\]
it results from (\ref{convCD}) that
\be
\label{convCDx}
E\{ \| \bm{x}_{CD}^{(k)} - \bm{x}_o \|^2 \} \le 
 \left( 1 - \frac{1}{\kappa_F^2( \bm{A} )} \right)^k
 \| \bm{A}^\dagger \|^2 \cdot \| \bm{b} - \bm{r}_o \|^2,
\ee

On the other side, \cite[Th.4.1]{ZoFr13} shows that the EK solution satisfies
\be
\label{convEK}
E\{ \| \bm{x}_{EK}^{(k)} - \bm{x}_o \|^2 \} \le
  \left( 1 - \frac{1}{\kappa_F^2( \bm{A} )} \right)^{\lfloor k/2 \rfloor} C,
\ee
where $C$ is a constant of no interest here.

Although the bounds (\ref{convCDx}) and (\ref{convEK}) are not necessarily tight,
they suggest that CD converges faster than EK, supporting the results of
Propositions \ref{prop:term} and \ref{prop:res}.
\QED
\end{remark}

%---------------------------------------------------------------------
\section{The general LS problem}

Reminding that all the previous discussion was for full-rank
overdetermined LS problems, let us look at the other cases.
Consider first underdetermined systems, but still full-rank.
In this case, the system $\bm{A}\bm{x} = \bm{b}$ is consistent
and it is well known that the standard Kaczmarz algorithm converges to the LS solution.
There is no need of residual approximation, since the residual is zero,
hence the CD part of EK is useless.

We pass now to the general LS problem, in which the matrix $\bm{A}$
has arbitrary rank, and for which the deterministic EK algorithm \cite{Popa95}
was originally intended.
The LS solution is that with minimum norm $\|\bm{x}\|_2$,
among those minimizing the residual $\| \bm{b} - \bm{A} \bm{x} \|_2$.
In this case, due to their specific projection operations,
CD can minimize the residual, but not find a solution with
minimum norm, while Kaczmarz can minimize the norm of the solution
(if properly initialized with $\bm{x}^{(0)}\in {\cal R}(\bm{A}^T)$),
but not that of the residual.
EK combines their strengths to find the LS solution.

We argue that, however, there are better ways to combine the two algorithms
than intertwining them as in EK.
We propose to run first CD for estimating the (nearly) optimal residual
$\bm{r}_{CD} \approx \bm{r}_o$
and only then Kaczmarz for finding the least norm solution of the
consistent system
\be
\bm{A}\bm{x} = \bm{b} - \bm{r}_{CD}.
\label{rsys}
\ee
We name CD+K this algorithm.
We note that the general idea of running CD before Kaczmarz is mentioned
in \cite{ZoFr13} (where CD is named "orthogonal projection");
however, the authors settle there for the specific form of EK and discuss
CD and K only separately.

\begin{remark} \rm
\label{rem:CDK}
In average, CD+K should need less Kaczmarz iterations than EK.
%Since the operations in EK can be perfectly split between CD and Kaczmarz,
%we count one iteration of CD and K together for one iteration of EK.
We cannot give a rigorous proof, but give below two heuristic
arguments supporting the above assertion.

{\em Argument 1.}
Since the CD operations are independent of the other operations in EK, it takes
the same number of iterations for CD and EK to satisfy the stopping criterion
(\ref{stopCD}).
Running Kaczmarz after CD has the advantage that it
works from the start on the (nearly) consistent system to be solved.
In EK, the Kaczmarz iterations are made for the system (\ref{Ksys}),
which is only an approximation of the final consistent system (\ref{rsys}).

More intuitively, while CD goes straightly to its target,
the Kaczmarz part of EK takes a detour.
Running first CD should be more efficient because
CD sets the final target, then K goes directly to it.
Both CD and K use their full power.
So, it is natural to expect less Kaczmarz iterations in CD+K than in EK.

{\em Argument 2.}
Let us take a look at the convergence speed.
In \cite[Th.3.4]{ZoFr13}, the Kaczmarz algorithm is shown to satisfy the relation
\be
\label{convK}
E\{ \| \bm{x}_{K}^{(k)} - \bm{x}_o \|^2 \} \le 
 \left( 1 - \frac{1}{\kappa_F^2( \bm{A} )} \right)^k
 \| \bm{x}_{K}^{(0)} - \bm{x}_o \|^2.
\ee
The constant bounding the convergence speed is the same as for CD, see (\ref{convCDx}).
So, the discussion from Remark \ref{rem:convCD-EK} applies also here,
suggesting that CD+K has better convergence speed than EK.
\QED
\end{remark}

\begin{remark} \rm
\label{rem:CDEKK}
One may argue that EK still has an advantage over CD+K: rough approximations of the
solution are earlier available. Indeed, in CD+K we have to wait for the
whole CD part before approximations of the solution are computed.
A possible fix is to recognize that between EK and CD+K there is a whole
family of algorithms, organized as follows.

First CD is run with a tolerance $\hat{\varepsilon}_{CD} > \varepsilon_{CD}$.
Then EK is run, initialized with the residual produced by CD,
until one of the stopping criteria (\ref{stopCD}) or (\ref{stopK}) is satisfied.
If (\ref{stopK}) is satisfied we stop, otherwise we run Kaczmarz on the now nearly
consistent system, initializing with the solution produced by EK, 
until (\ref{stopK}) is met.
We name CD+EK+K this algorithm.
Again, we expect it to be faster than EK, the arguments being similar
to those in Remark \ref{rem:CDK}.
If $\hat{\varepsilon}_{CD}$ is large, only few CD iterations are made,
hence approximations of the solution are quickly available.
\QED
\end{remark}

\begin{remark} \rm
Running Kaczmarz after CD has a nice alternative interpretation.
In this context, the Kaczmarz iteration has the form
(see step 8 of Algorithm 1)
\be
\bm{x}_{K}^{(k+1)} = \bm{x}_{K}^{(k)} +
     \frac{\langle \bm{b} - \bm{r}_{CD}, \bm{e}_{i_k} \rangle
              - \langle \bm{x}_{K}^{(k)},\bm{A}^{(i_k)}\rangle}
          {\|\bm{A}^{(i_k)}\|_2^2} \bm{A}^{(i_k)}
 = \bm{x}_{K}^{(k)} +
     \frac{\langle \bm{x}_{CD} - \bm{x}_{K}^{(k)},\bm{A}^{(i_k)}\rangle}
          {\|\bm{A}^{(i_k)}\|_2^2} \bm{A}^{(i_k)},
\label{KafterCD}
\ee
where we have used the equality
$\langle \bm{b} - \bm{r}_{CD}, \bm{e}_{i_k} \rangle
= \langle \bm{A} \bm{x}_{CD}, \bm{e}_{i_k} \rangle = 
\langle \bm{x}_{CD}, \bm{A}^{(i_k)} \rangle$.
Denoting $\bm{q}_{K}^{(k)} = \bm{x}_{CD} - \bm{x}_{K}^{(k)}$, it results from
(\ref{KafterCD}) that
\be
\bm{q}_{K}^{(k+1)} = \bm{q}_{K}^{(k)} -
     \frac{\langle \bm{q}_{K}^{(k)},\bm{A}^{(i_k)}\rangle}
          {\|\bm{A}^{(i_k)}\|_2^2} \bm{A}^{(i_k)}.
\label{Kdual}
\ee
This is a projection operation on the orthogonal complement of the $i_k$-th row
of $\bm{A}$, dual to the CD operation of projecting the residual on the
orthogonal complement of column $j_k$ (step 5 of Algorithm 1).
Hence, $\bm{q}_{K}^{(k)}$ tends to the component of $\bm{q}_{K}^{(0)}$
that is orthogonal on ${\cal R}(\bm{A}^T)$.
Initializing with $\bm{q}_{K}^{(0)} = \bm{x}_{CD}$, which corresponds
to the natural initialization $\bm{x}_{K}^{(0)} = 0$, the iteration
(\ref{Kdual}) converges to $\bm{q}_K$ satisfying
$\bm{x}_K + \bm{q}_K = \bm{x}_{CD}$, with $\bm{q}_K \perp {\cal R}(\bm{A}^T)$
and hence $\bm{x}_K  \in {\cal R}(\bm{A}^T)$.
Since $\bm{A} \bm{q}_K = 0$, the Kaczmarz solution satisfies the system (\ref{rsys}).
This means that $\bm{x}_K = \bm{x}_o$, since the LS solution is the unique
vector from ${\cal R}(\bm{A}^T)$ satisfying (\ref{rsys}).

So, the iteration (\ref{Kdual}) starts with $\bm{x}_{CD}$, for which (\ref{rsys})
already holds, and computes its projection onto the orthogonal complement
of ${\cal R}(\bm{A}^T)$.
Thus, it allows the computation of the projection of $\bm{x}_{CD}$ onto
${\cal R}(\bm{A}^T)$, which is the LS solution.

Of course, using (\ref{Kdual}) instead of (\ref{KafterCD}) gives no computational
advantage, but gives a dual view to the convergence of the Kaczmarz iterations.
\QED
\end{remark}

%---------------------------------------------------------------------
\section{Numerical results}

We have implemented Algorithm 1 in Matlab and report the performance
of CD and EK only in terms of number of iterations, reminding however
that at similar number of iterations CD is still faster.
The algorithm has been run for a sufficiently high number of iterations,
without any stopping criterion.
For overdetermined LS problems, we have considered matrices
belonging to two classes where EK was shown
in \cite{ZoFr13} to have better performance than other algorithms:
(i) dense well-conditioned matrices, generated with {\tt randn},
and (ii) sparse random matrices with density $0.25$, generated
with {\tt sprandn}.

\begin{figure}
  \centering
  \includegraphics[scale=0.4]{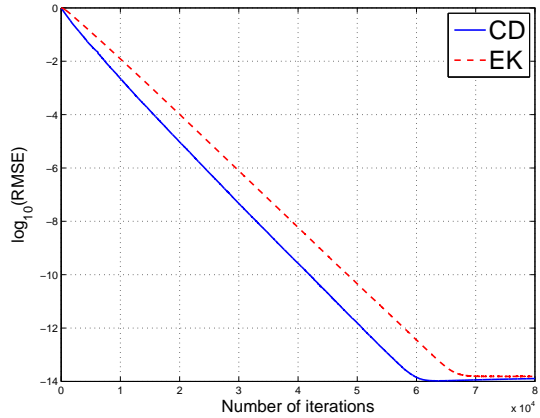}\\
  \caption{RMSE for dense matrices, $m=2000$, $n=500$.}
  \label{fig:1}
\end{figure}

\begin{figure}
  \centering
  \includegraphics[scale=0.4]{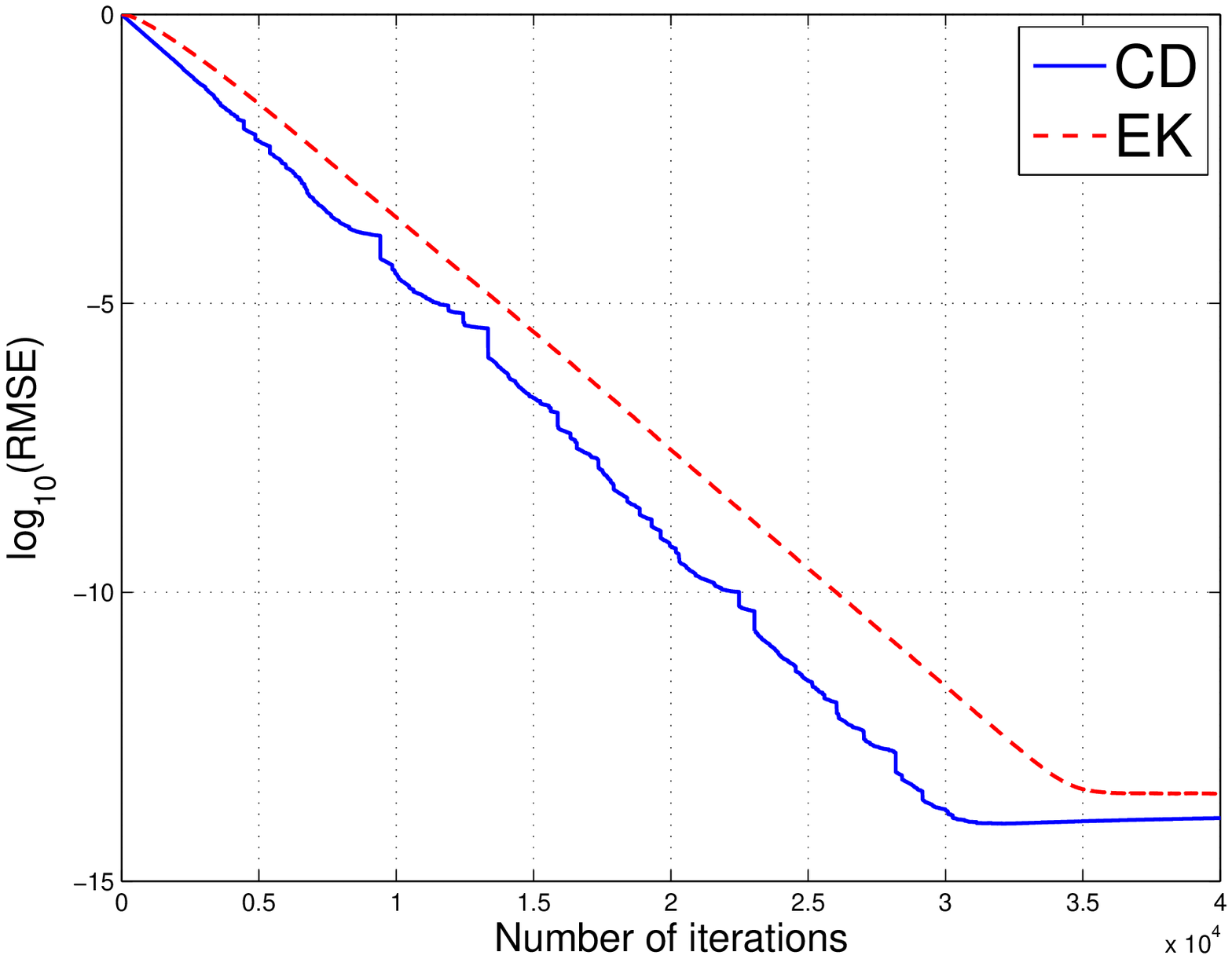}\\
  \caption{RMSE for dense matrices, $m=10000$, $n=500$.}
  \label{fig:3}
\end{figure}

\begin{figure}
  \centering
  \includegraphics[scale=0.4]{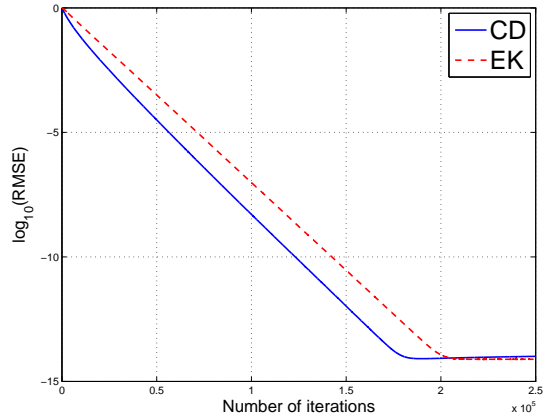}\\
  \caption{RMSE for sparse matrices, $m=2000$, $n=800$.}
  \label{fig:2}
\end{figure}

We report the normalized RMSE
$\sqrt{ E\left( \| \bm{x}^{(k)} - \bm{x}_o \|^2 / \| \bm{x}_o \|^2 \right) }$,
obtained by averaging over 100 matrices from the same class.
Figures \ref{fig:1} and \ref{fig:3} show the RMSE for dense matrices with
the same number of rows, $n=500$, but different number of columns:
$m=2000$ and $m=10000$, respectively.
Figure \ref{fig:2} shows the RMSE for sparse matrices, $n=800$, $m=2000$.
In all cases, the faster convergence of CD is clear.
When the system is very overdetermined, CD has a jagged convergence,
alternating many small advances with few large ones, but is still faster.
For other matrix sizes, the results are similar.

%--------------------------
%--------------------------
%--------------------------
%--------------------------
%--------------------------

To illustrate the behavior of CD+K and CD+EK+K, we have generated random matrices
$\bm{A} \in \R^{m \times n}$, computed their SVD
$\bm{A} = \bm{U} \bm{\Sigma} \bm{V}^T$, kept only the $r$ largest singular
values in $\bm{\Sigma}$ while setting the others to zero,
and recomputed $\bm{A}$ using the same above relation.
We have implemented CD+EK+K in a simple manner that ensures a
simple evaluation: first CD is run for $N$ iterations, then EK for a
number of iterations (that may vary depending on matrix size and rank),
then K is run for $N$ iterations.
Note that CD and K run the same number of iterations, like in EK. 

\begin{figure}
  \centering
  \includegraphics[scale=0.4]{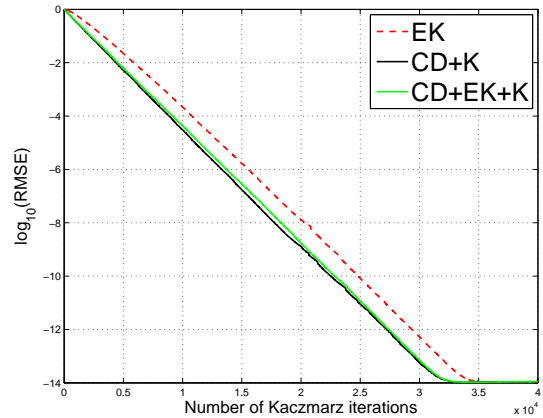}\\
  \caption{RMSE for dense matrices, $m=500$, $n=2000$, $N=2000$.}
  \label{fig:g1}
\end{figure}

We give only one typical sample of result, for an underdetermined problem with
$m=500$, $n=2000$, $r=400$.
Figure \ref{fig:g1} shows the RMSE of the solution approximations
as a function of the number of Kaczmarz iterations; the number of
CD iteration is the same for all methods, but their position in time is different.
It is visible that both CD+K and CD+EK+K need less Kaczmarz iterations to converge.
Although $N = 2000$ for CD+EK+K, which is a relatively small value,
the performance is about the same as for CD+K; for larger values, like
$N=5000$, the curves for CD+K and CD+EK+K are nearly identical;
for smaller values, the curve of CD+EK+K approaches that of EK.

All the presented curves show that the advantage of CD, CD+K or CD+EK+K
over EK is built especially in the first iterations.
This corresponds well with the fact that in those iterations the
approximation of the residual is still very poor in EK, and hence the
Kaczmarz iterations do not have a good target, as explained in
Remark \ref{rem:CDvsEK} and in argument 1 of Remark \ref{rem:CDK}.

%--------------------------
%--------------------------

%---------------------------------------------------------------------
\section{Conclusions}

The computational conclusion of all the facts presented in this note
is the recommendation to replace EK with one of the following algorithms: 
\bi
 \item CD, for full-rank overdetermined LS problems;
 \item CD+K or CD+EK+K, for rank-deficient or unknown rank problems.
\ei
(As already known, Kaczmarz replaces EK for full-rank underdetermined problems.)

In particular, for the full-rank overdetermined problem,
CD is always preferable to EK, due to the following reasons:

\begin{itemize}
 \item In average, CD converges in less iterations than EK.
 \item CD needs less operations per iteration.
 \item CD uses only the columns of the matrix $\bf{A}$,
while EK uses both columns and rows.
\end{itemize}

The conclusions apply equally to randomized EK \cite{ZoFr13}
and the less efficient deterministic version \cite{Popa95}.

%---------------------------------------------------------------------
\section*{Acknowledgment}

The author is indebted to Liang Dai for stimulating discussions on the
Kaczmarz algorithm and for pointing out \cite{ZoFr13}.
%and to the reviewers for their constructive comments.

%---------------------------------------------------------------------
%\bibliographystyle{plain}
%\bibliography{../../biblio/random,../../biblio/mybib}

%---------------------------------------------------------------------

\end{document}